# ON THE DESIGN-CONSISTENCY PROPERTY OF HIERARCHICAL BAYES ESTIMATORS IN FINITE POPULATION SAMPLING

By P. Lahiri[1] and Kanchan Mukherjee

*University of Maryland, College Park and University of Liverpool*

We obtain a limit of a hierarchical Bayes estimator of a finite population mean when the sample size is large. The limit is in the sense of ordinary calculus, where the sample observations are treated as fixed quantities. Our result suggests a simple way to correct the hierarchical Bayes estimator to achieve design-consistency, a well-known property in the traditional randomization approach to finite population sampling. We also suggest three different measures of uncertainty of our proposed estimator.

**1. Introduction.** Ericson [4] put forward a subjective Bayesian approach in finite population sampling. The subjective Bayes and the more general hierarchical Bayes estimators have been found to be effective in combining information from a variety of sources in conjunction with the sample survey data on the main variable of interest; see [10]. Specific applications include small area estimation, estimation from longitudinal surveys, and so on.

In Section 2, we obtain the mathematical limit of a general version of Ericson's subjective Bayes estimator of a finite population mean when the sample size is large. By mathematical limit, we mean the limit in the sense of ordinary calculus, where observations are treated as fixed real numbers. We show that the Bayes estimator converges to a quantity which is free of any hyperparameters, but may depend on the sample observations, irrespective of the model used to derive the Bayes estimator. This pure mathematical result can be used to examine the design-consistency property of the subjective Bayes estimator. Design-consistency is a desirable property in the

Received March 2003; revised May 2006.
[1]Supported in part by a grant from the Center of Excellence in Health Statistics of the University of Michigan, Ann Arbor, and the Gallup Organization.

randomization approach to finite population sampling. Kott [13] advocated the use of a model-based estimator of a finite population mean which is also design-consistent. For a formal definition of design-consistency, see [22], page 18. We find that the subjective Bayes estimator is, in general, not design-consistent. In other words, the subjective Bayes estimator does not converge to the true finite population mean as the sample size becomes large. Here the convergence is defined with respect to the sampling design, observations for all units of the finite population being treated as fixed non-random quantities. The result, in turn, suggests a simple correction to the Bayes estimator to achieve design-consistency.

In Section 3 we consider the case where the finite population is divided into several strata and a sample is available from each stratum. Model-based estimation (both empirical and hierarchical Bayes) of a stratum mean has been considered by several authors; see Ghosh and Meeden [9], Ghosh and Lahiri [6], Battese, Harter and Fuller [2], Datta and Ghosh [3], Prasad and Rao [19], Ghosh and Lahiri [8], Arora, Lahiri and Mukherjee [1] and Jiang, Lahiri and Wan [12], among others. Prasad and Rao [20] observed that these model-based estimators are typically not design-consistent unless the sampling design is self-weighting within each stratum, and they proposed a pseudo-EBLUP design-consistent estimator for the normal case; see [23] for a pseudo-hierarchical Bayes version of Prasad and Rao [20]. However, it is not clear how one can extend their approach to a general nonnormal case, since the distribution of a linear combination of observations used in their paper is not always analytically tractable for nonnormal cases. In Section 3 we provide the stratified sampling extension of the result of Section 2 for a general hierarchical Bayes estimator of the finite population mean when the hyperparameters are assumed to be unknown. Again, our mathematical limit result lends itself to a simple correction of the hierarchical Bayes estimator to achieve the design-consistency property.

For small area problems with binary data, Folsom, Shah and Vaish [5] proposed a method which appears to work well in their simulation work in terms of design-consistency. However, a formal proof is needed to claim design-consistency of their estimator. In this connection, we also mention the work of Malec, Davis and Cao [16], who considered sample selection adjustment to their model-based procedure. See [21] and [11] for further discussions of design-consistency in small area estimation. The simple and general approach taken in this paper merely involves finding the mathematical limit of a hierarchical Bayes estimator.

Measuring uncertainty of the proposed design-consistent model-based estimator is an important problem. In Section 4 we first present two possible ways to measure the uncertainty of the proposed design-consistent model-based estimator. We note their merits and disadvantages and propose a third measure that is just a simple average of the two. The results from a simulation study are presented in Section 5. In our simulation experiment the proposed measures perform reasonably well and the third measure appears to be a compromise between the first two. The technical derivations of the results of Sections 2 and 3 are relegated to the Appendix.

**2. The Bayes estimator and its limiting behavior.** Let $Y_i$ be the value of the characteristic of interest for the $i$th unit of a finite population ($i = 1, \ldots, N$). We assume that the finite population size $N$ is known. In our subjective Bayesian formulation, we assume that $Y_i$, $i = 1, \ldots, N$, are independent realizations from a superpopulation belonging to the family of densities

$$(2.1) \qquad f(y|\theta, \phi) = \exp[\phi^{-1}\{y\theta - \psi(\theta)\} + \rho(y, \phi)],$$

where $\phi > 0$. When $\phi$ is known, the generalized linear models characterized by the above class of densities include the exponential family of distributions (in particular, Gaussian, Bernoulli and Poisson distributions). Thus, it includes continuous data as well as both the categorical and count data. When $\phi$ is unknown, it includes distributions which are not in the exponential family. Note that the function $\psi$ is given by

$$\psi(\theta) = \phi \log \int \exp[\phi^{-1} y\theta + \rho(y, \phi)] \, dy$$

and the mean of $Y$ is $\psi'(\theta)$; see [18], page 28. Here and throughout the paper, $\alpha'$ denotes the derivative of an arbitrary function $\alpha$.

We assume that the superpopulation parameter $\theta$ has the prior distribution

$$(2.2) \qquad h(\theta) = \beta + u,$$

where $h$, called the *link function*, is a strictly increasing function of $\theta$, $\beta \in \mathbb{R}$ is a location parameter and $u \sim N(0, r^{-1})$. Similar Bayesian models were used by Ghosh and Meeden ([10], page 269).

A sample of fixed size $n$ is drawn from the finite population. Let $\{p(s)\}$ denote the sampling design. Note that $p(s)$ is the probability of drawing a particular sample $s$ of size $n$ from the universe of all possible samples $S$ of size $n$. Thus, $p(s) \geq 0$ for all $s \in S$ and $\sum_{s \in S} p(s) = 1$. Let $\pi_i = P(s \ni i)$ be the first-order inclusion probability of unit $i$, that is, the probability of including the population unit $i$ in the sample ($i = 1, \ldots, N$). Let $y_s = (Y_i, i \in s)$ be the vector of observations in the sample.

The standard Bayesian approach to finite population sampling recognizes the importance of a good sampling design for selecting the sample, but once the sample is selected, the approach does not use the inclusion probabilities in the estimation. Under the Bayesian model (2.1) and (2.2) and squared error loss, the Bayes estimator of the finite population mean, $\overline{Y} = \sum_{i=1}^{N} Y_i/N$, is given by

$$\widehat{\overline{Y}}^{\mathrm{B}} = E[\overline{Y}|y_s] = f_n \bar{y}_s + (1-f_n)E[\psi'(\theta)|y_s], \tag{2.3}$$

where $f_n = n/N$ and $\bar{y}_s = n^{-1}\sum_{i\in s} Y_i$ is the sample mean.

In order to understand the relationship between the Bayes estimator $\widehat{\overline{Y}}^{\mathrm{B}}$ and $\bar{y}_s$ for a large sample, we find the mathematical limit of $E[\psi'(\theta)|y_s]$ as $n \to \infty$ as a *function* of $\bar{y}_s$. Hence, in the following $\bar{y}_s$ is treated as a nonrandom argument of $E[\psi'(\theta)|y_s]$. Consequently, let $T_n$ be a random variable with density (as a function of $t$)

$$\exp[-n\phi^{-1}\psi(t)]\exp[\phi^{-1}n\bar{y}_s t] \Big/ \int \exp[-n\phi^{-1}\psi(\tau)]\exp[\phi^{-1}n\bar{y}_s \tau]\,d\tau. \tag{2.4}$$

This leads to the following lemma.

LEMMA 2.1. *In addition to* (2.1) *and* (2.2), *assume the following regularity conditions:*

(R.1) *The functions* $\psi'(\cdot)\exp[-(r/2)\{h(\cdot)-\beta\}^2]h'(\cdot)$ *and* $\exp[-(r/2)\{h(\cdot)-\beta\}^2]h'(\cdot)$ *are bounded.*

(R.2) *The sequence of random variables* $\{T_n\}$ *converges in probability with respect to* (2.4) *to* $T$ *as* $n \to \infty$.

Then $E[\psi'(\theta)|y_s]$ *converges to* $C := \psi'(T)$.

REMARK 2.1. In general, $T$ is a function of $\bar{y}_s$. This reveals the large-sample behavior of $E[\psi'(\theta)|y_s]$ as a function of $\bar{y}_s$. In the Appendix, we verify that assumptions (R.1) and (R.2) are satisfied for three well-known distributions involving Gaussian, Bernoulli and Poisson distributions, where $C = \bar{y}_s$.

The following theorem is a simple consequence of Lemma 2.1 and (2.3).

THEOREM 2.1. *In addition to the conditions of Lemma* 2.1, *suppose that* $f_n \to f$ *for some* $0 < f < 1$. *Then*

$$\widehat{\overline{Y}}^{\mathrm{B}} \to f\bar{y}_s + (1-f)C$$

*as* $n \to \infty$.

In the Appendix, we note that the convergence proof of Theorem 2.1 does not use any assumptions regarding the sampling design. In general, for large $n$, $\widehat{\overline{Y}}^{\mathrm{B}}$ is not design-consistent, except possibly for a self-weighting sampling design. Let $\bar{y}_w$ be any design-consistent estimator of $\overline{Y}$. For example, we can choose $\bar{y}_w$ to be the well-known Hansen–Hurwitz or Horvitz–Thompson estimator. Then, using Theorem 2.1, we can obtain the following design-consistent estimator of $\overline{Y}$ based on the Bayes estimator:

$$\widehat{\overline{Y}} = \widehat{\overline{Y}}^{\mathrm{B}} - \{f\bar{y}_s + (1-f)C - \bar{y}_w\}.$$

For a self-weighing design, $\bar{y}_w = \bar{y}_s$, and so the Bayes estimator is indeed design-consistent for the three examples given in the Appendix.

**3. Hierarchical Bayes estimator and its limiting behavior.** In this section we extend the results of Section 2 and consider two important cases: *case* (i) $\phi$ is known, but $\beta$ and $r$ are unknown, and *case* (ii) all of the hyperparameters are unknown. To treat these two cases, we need further information on the finite population that allows estimation of unknown hyperparameters through an appropriate hierarchical Bayes method. To this end, consider a finite population divided into $m$ strata. Let $Y_{ij}$ denote the value of the $j$th observation in the $i$th stratum $(i=1,\ldots,m; j=1,\ldots,N_i)$. We consider the estimation of a particular stratum mean. Without loss of generality, we consider estimation of the $m$th stratum mean given by $\overline{Y}_m = N_m^{-1}\sum_{j=1}^{N_m} Y_{mj}$, where $N_m$ is the known population size for the $m$th stratum. We assume that the sample $y_s$ consists of $n_i$ observations from the $i$th stratum, $1 \leq i \leq m$.

In addition to (2.1) and (2.2), for each $Y_{ij}$ $(i=1,\ldots,m; j=1,\ldots,N_i)$, we assume that $\beta$ and $r$ are independent, with $\beta \sim U(-\infty,\infty)$, an improper uniform distribution over the real line, and $r \sim G(a,b)$, a gamma distribution with density proportional to $r^{(b/2)-1}\exp[-(a/2)r]$.

Writing $\boldsymbol{\theta} = [\theta_1,\ldots,\theta_m]'$, we obtain the conditional density

$$\pi(\boldsymbol{\theta},\boldsymbol{\beta},r,\phi|y_s) \propto \prod_{i=1}^{m}\prod_{j=1}^{n_i}\exp[\phi^{-1}\{y_{ij}\theta_i - \psi(\theta_i)\} + \rho(y_{ij},\phi)]$$

(3.1)
$$\times \prod_{i=1}^{m}\{\exp[-(r/2)\{h(\theta_i) - \beta\}^2]h'(\theta_i)r^{1/2}\}$$

$$\times \exp[-ar/2]r^{(b/2)-1}.$$

First, consider *case* (i). In this case, we show that (see the Appendix)

(3.2)
$$E[\psi'(\theta_m)|y_s]$$
$$= E_{T_{n_m}}[\psi'(T_{n_m})h'(T_{n_m})g(T_{n_m})]/E_{T_{n_m}}[h'(T_{n_m})g(T_{n_m})],$$

where the expectation in (3.2) is taken with respect to a random variable $T_{n_m}$ with density (as a function of $t$)

(3.3) $\exp[\phi^{-1}n_m\{\bar{y}_m t - \psi(t)\}]\Big/\int \exp[\phi^{-1}n_m\{\bar{y}_m\tau - \psi(\tau)\}]\,d\tau,$

and where $g$ is a function defined in (A.7).

This leads to the following lemma, which, in turn, proves the convergence of the hierarchical Bayes estimator to $\bar{y}_m$.

LEMMA 3.1. *Assume the following regularity conditions.*

(R.3) *The functions $\psi'(\cdot)h'(\cdot)g(\cdot)$ and $h'(\cdot)g(\cdot)$ are bounded.*

(R.4) *The sequence of random variables $\{T_{n_m}\}$ converges in probability [with respect to (3.3)] to $T$ as $n_m \to \infty$.*

*Then $E[\psi'(\theta_m)|y_s]$ converges to $C_m = \psi'(T)$.*

REMARK 3.1. As before, we can check that the regularity conditions (R.3) and (R.4) are satisfied for all three examples in the Appendix and that in each case $C_m = \bar{y}_m$. Verification of (R.4) is similar to that of (R.2). Verification of (R.3) may be facilitated by noting that

$$g(\theta_m) \leq \int \prod_{i=1}^{m-1} \left\{ h'(\theta_i) \prod_{j=1}^{n_i} \exp[\phi^{-1}\{y_{ij}\theta_i - \psi(\theta_i)\}] \right\}$$
$$\times \left[ a + \sum_{i=1}^{m} h^2(\theta_i) \right]^{-(b+m-1)/2} \prod_{i=1}^{m-1} d\theta_i.$$

For example, in the case of normal distribution,

$$g(\theta_m) \leq [a + \theta_m^2]^{-(b+m-1)/2} \int \prod_{i=1}^{m-1} \prod_{j=1}^{n_i} \exp[\sigma^{-2}\{y_{ij}\theta_i - \theta_i^2/2\}] \prod_{i=1}^{m-1} d\theta_i < \infty.$$

The following theorem is an immediate consequence of Lemma 3.1. The theorem suggests a simple adjustment to the hierarchical Bayes estimator to achieve the design-consistency property.

THEOREM 3.1. *In addition to the regularity conditions of Lemma 3.1, assume that $f_{n_m} := n_m/N_m \to f_m$ for some $0 < f_m < 1$. Let $\widehat{\overline{Y}}_m^{HB}$ be the hierarchical Bayes estimator of $\overline{Y}_m$. Then*

$$\widehat{\overline{Y}}_m^{HB} \to f_m \bar{y}_m + (1 - f_m) C_m$$

*as $n_m \to \infty$.*

Let $\bar{y}_{mw}$ be any design-consistent estimator of $\overline{Y}_m$. Then $\widehat{\widehat{Y}}_m = \widehat{\overline{Y}}_m^{HB} - \{f_m \bar{y}_m + (1 - f_m)C_m - \bar{y}_{mw}\}$ is design-consistent, based on the hierarchical Bayes estimator.

Now consider *case* (ii). Suppose $\phi$, $\beta$ and $r$ are mutually independent. Furthermore, let $v := 1/\phi \sim U(0, \infty), \beta \sim U(-\infty, \infty)$ and $r \sim G(a, b)$.

In this case, we will consider the normal distribution only, since $\phi$ is known for the binomial and Poisson examples. Since $h(\theta) = \theta$, integrating (3.1) with respect to $\beta$, $r$ and $v$ (in that order), we obtain

(3.4)
$$\pi(\boldsymbol{\theta}|y_s) \propto \left\{\sum_{i=1}^{m}\sum_{j=1}^{n_i}(y_{ij} - \theta_i)^2\right\}^{-(n_T/2+1)} \left[a + \sum_{i=1}^{m}\{\theta_i - \bar{\theta}\}^2\right]^{-(b+m-1)/2},$$

where $n_T = \sum_{i=1}^{m} n_i$. In the Appendix, we show that

(3.5) $$E[\psi'(\theta_m)|y_s] \to \bar{y}_m.$$

**4. Uncertainty measures.** A conventional measure of uncertainty of the hierarchical Bayes estimator $\widehat{\overline{Y}}_m^{\text{HB}}$ is simply the posterior variance $V(\overline{Y}_m|y_s)$. Noting that $V(\overline{Y}_m|y_s) = E[(\widehat{\overline{Y}}_m^{\text{HB}} - \overline{Y}_m)^2|y_s]$, we may define a measure of uncertainty (MU) of $\widehat{\overline{Y}}_m$ as

$$MU_m^{(1)} = E[(\widehat{\overline{Y}}_m - \overline{Y}_m)^2|y_s] = V(\overline{Y}_m|y_s) + [\widehat{\overline{Y}}_m - \widehat{\overline{Y}}_m^{\text{HB}}]^2.$$

Thus, in order to achieve design-consistency, we increase this measure by $[\widehat{\overline{Y}}_m - \widehat{\overline{Y}}_m^{\text{HB}}]^2$. However, this apparent increase may be misleading, since this will only happen if the assumed hierarchical model holds for all units of the finite population, an assumption hard to justify for the unobserved units of the finite population based on the observed units in the sample.

We now propose an alternative measure of uncertainty following the work of Prasad and Rao [20], who considered a design-consistent pseudo-EBLUP for a nested error regression model. Following their approach, we define the mean squared error of $\widehat{\overline{Y}}_m$ as

$$MSE(\widehat{\overline{Y}}_m) = E[\widehat{\overline{Y}}_m - \eta_m]^2,$$

where $\eta_m = E(\overline{Y}_m|\theta_m)$, the expectation being taken over the first two levels of the hierarchical model. Unlike the previous approach, this approach does not require extensive model assumptions regarding the unobserved units of the finite population, except for the mild assumption of the existence of a random effect $\theta_m$. This is certainly an advantage of this approach over the previous approach.

Let $\hat{\eta}_m^{\text{B}} = \hat{\eta}_m^{\text{B}}(\phi)$ be the Bayes estimator of $\eta_m$. Note that

$$MSE(\widehat{\overline{Y}}_m) = h_{1m}(\phi) + h_{2m}(\phi),$$

where $h_{1m}(\phi) = E[\hat{\eta}_m^{\text{B}} - \eta_m]^2$ and $h_{2m}(\phi) = E[\widehat{\overline{Y}}_m - \hat{\eta}_m^{\text{B}}]^2$. It is possible to write down an explicit expression for $h_{1m}(\phi)$, although it may not be in closed form.

Let $\hat{\phi}$ be any commonly used consistent estimator of $\phi$. For example, in a mixed linear normal model, we can consider the residual maximum likelihood estimators (REML) for the variance components and weighted least squares with estimated variance components for the regression coefficients. We can estimate $h_{1m}(\phi)$ by $h_{1m}(\hat{\phi})$ and $h_{2m}(\phi)$ by $\hat{h}_{2m}(\hat{\phi}) = [\widehat{\overline{Y}}_m - \hat{\eta}_m^{\text{EB}}]^2$, where $\hat{\eta}_m^{\text{EB}} = \hat{\eta}_m^{\text{B}}(\hat{\phi})$, an empirical Bayes estimator of $\eta_m$. We propose the following as the second measure of uncertainty of $\widehat{\overline{Y}}_m$:

$$MU_m^{(2)} = h_{1m}(\hat{\phi}) + \hat{h}_{2m}(\hat{\phi}).$$

Note that $MU_m^{(2)}$ does not incorporate the variability due to the estimation of $\phi$. On the other hand, $MU_m^{(1)}$ incorporates all sources of variability, but may be sensitive to the specification of the prior distribution on $\phi$. Thus, as a compromise, we propose the following measure of uncertainty:

$$MU_m^{(3)} = \frac{MU_m^{(1)} + MU_m^{(2)}}{2}.$$

**5. A Monte Carlo simulation.** In this section, we conduct a Monte Carlo simulation study to compare the performances of $MU^{(1)}$, $MU^{(2)}$ and $MU^{(3)}$, proposed in Section 4. In particular, we study the performances of these measures in estimating the design-based mean squared error defined as $MSE_d(\widehat{\overline{Y}}_m) = E_d(\widehat{\overline{Y}}_m - \overline{Y}_m)^2$, where $E_d$ denotes an expectation with respect to the sampling design.

We consider 100 finite populations, each of size 60. The main variable $Y$ is generated for each unit of the finite populations using a nested error model $Y_{ij} = \mu + v_i + e_{ij}$, where $\mu$ is the fixed effect, $v_i \sim N(0, \sigma_v^2)$ are the random effects and $e_{ij} \sim N(0, \sigma_e^2)$, with the $\{v_i\}$'s and the pure errors $\{e_{ij}\}$'s assumed to be independent, $1 \le i \le 100$, $1 \le j \le N_i = 60$. We set $\mu = 50$ and $\sigma_v = 1$ and consider two different values of $\sigma_e$, namely $\sigma_e = 1$ and 2.

We draw a sample of size $n$ from each finite population using a probability proportional to size with replacement (PPSWR) sampling design. We consider three different choices of $n$, namely 10, 20 and 30. The size measure is generated for each unit of the finite populations using an exponential distribution with mean 1.

Flat priors on the hyperparameters are used to obtain the posterior mean and the posterior variance needed to compute $MU^{(1)}$. We use PROC MIXED in SAS to generate 1050 observations from the posterior distributions, but only the last 1000 observations are retained for approximating the posterior means and variances. Essentially, PROC MIXED uses a Markov chain Monte Carlo (MCMC) technique. PROC IML is used to obtain the required posterior means and posterior variances. For the second measure $MU^{(2)}$, $\phi$ is estimated by the residual maximum likelihood (REML) method.

Table 1
*ARB and ARRMSE of the three different measures of uncertainty*

| $n$ | ARB | | | ARRMSE | | |
|---|---|---|---|---|---|---|
| | $ARB^{(1)}$ | $ARB^{(2)}$ | $ARB^{(3)}$ | $ARRMSE^{(1)}$ | $ARRMSE^{(2)}$ | $ARRMSE^{(3)}$ |
| $\sigma_e = 1$ | | | | | | |
| 10 | 0.001 | −0.051 | −0.025 | 1.132 | 0.523 | 0.743 |
| 20 | 0.011 | −0.004 | 0.004 | 1.340 | 0.673 | 0.891 |
| 30 | 0.017 | 0.010 | 0.013 | 1.515 | 0.761 | 1.003 |
| $\sigma_e = 2$ | | | | | | |
| 10 | 0.052 | −0.109 | −0.029 | 1.266 | 0.726 | 0.935 |
| 20 | 0.045 | −0.060 | −0.007 | 1.471 | 0.821 | 1.051 |
| 30 | 0.044 | −0.022 | 0.011 | 1.566 | 0.883 | 1.101 |

The relative bias (RB) and relative root mean square error (RRMSE) for the $i$-population ($1 \leq i \leq 100$) using the $k$th measure of uncertainty ($1 \leq k \leq 3$) are defined as

$$RB_i^{(k)} = \frac{E_d[MU_i^{(k)}] - MSE_d(\widehat{\overline{Y}}_i)}{MSE_d(\widehat{\overline{Y}}_i)}$$

and

$$RRMSE_i^{(k)} = \frac{\sqrt{E_d[MU_i^{(k)} - MSE_d(\widehat{\overline{Y}}_i)]^2}}{MSE_d(\widehat{\overline{Y}}_i)},$$

respectively. Table 1 reports the average RB $\{ARB^{(k)}; 1 \leq k \leq 3\}$ and average RRMSE $\{ARRMSE^{(k)}; 1 \leq k \leq 3\}$ for different combinations of $(n, \sigma_e)$, the average being taken over all 100 finite populations. In terms of ARB, it appears that the measure $MU^{(1)}$ has a slight tendency to overestimate the design-based MSE, whereas $MU^{(2)}$ has a slight tendency to underestimate. This is probably due to the fact that $MU^{(1)}$ attempts to incorporate all sources of variability, while $MU^{(2)}$ does not incorporate the variability in estimating $\phi$. It is interesting to note that for all the measures, ARRMSE increases with the increase of $n$. This behavior can be explained by the pattern of the inclusion probabilities induced by our sampling design; see [15] for details. In terms of the ARRMSE, $MU^{(2)}$ is better than both $MU^{(1)}$ and $MU^{(3)}$. Overall, the measure $MU^{(3)}$ is a compromise between $MU^{(1)}$ and $MU^{(2)}$.

**6. Concluding remarks.** In this paper, we examine a useful asymptotic behavior of the hierarchical Bayes estimator of a finite population mean. This leads to a simple method for constructing a design-consistent model-based estimator of a finite population mean. The method is general, in that it can be easily applied to both normal and nonnormal cases and is applicable to any complex weighting scheme. We have also addressed the important issue of measuring uncertainty of our proposed estimator. The simulation study reveals that our second measure suffers from a slight downward design-based bias. In the future, the Taylor series or a parametric bootstrap method as in [14] may be considered in an effort to reduce the bias. It is conceivable that our method extends beyond the stratified sampling design, for example, the multi-stage sampling design considered in [7] and [17], but this needs further research.

## APPENDIX

**Proof of Lemma 2.1.** First we compute $E[\psi'(\theta)|y_s]$. Using (2.1), (2.2) and the fact that if a random variable $h(\theta)$ has density $d(\cdot)$, then $\theta$ has density $d\{h(\theta)\}h'(\theta)$, we obtain the conditional density

$$
\begin{aligned}
\pi(\theta|y_s) &\propto \exp[-(r/2)\{h(\theta) - \beta\}^2]h'(\theta) \\
&\quad \times \exp[-n\phi^{-1}\psi(\theta)] \times \exp[\phi^{-1}n\bar{y}_s\theta].
\end{aligned}
\tag{A.1}
$$

This yields

$$
\begin{aligned}
&E[\psi'(\theta)|y_s] \\
&= \left\{\int \psi'(t)\exp[-(r/2)\{h(t) - \beta\}^2]h'(t)\exp[-n\phi^{-1}\psi(t)]\exp[\phi^{-1}n\bar{y}_s t]\,dt\right\} \\
&\quad \div \left\{\int \exp[-(r/2)\{h(t) - \beta\}^2]h'(t)\exp[-n\phi^{-1}\psi(t)]\exp[\phi^{-1}n\bar{y}_s t]\,dt\right\},
\end{aligned}
$$

which is a function of $n$ and $\bar{y}_s$.

Therefore,

$$
\begin{aligned}
E[\psi'(\theta)|y_s] &= \{E_{T_n}[\psi'(T_n)\exp[-(r/2)\{h(T_n) - \beta\}^2]h'(T_n)]\} \\
&\quad \div \{E_{T_n}[\exp[-(r/2)\{h(T_n) - \beta\}^2]h'(T_n)]\},
\end{aligned}
\tag{A.2}
$$

where the expectation in (A.2) is taken with respect to a random variable $T_n$ with density (2.4). The proof now follows by using the bounded convergence theorem on both the numerator and the denominator of (A.2).

*Verification of the conditions of Lemma* 2.1. When $\phi$ is known in (2.1), we verify the conditions of Lemma 2.1 for the Gaussian, Bernoulli and Poisson distributions. Verification of (R.1) is trivial.

EXAMPLE 1. Suppose that conditional on $\theta$, $Y$ is normal with mean $\theta$ and variance $\sigma^2$. Then equation (2.4) becomes

$$\exp[-nt^2/(2\sigma^2)]\exp[tn\bar{y}_s/\sigma^2] \Big/ \int \exp[-n\tau^2/(2\sigma^2)]\exp[\tau n\bar{y}_s/\sigma^2]\,d\tau,$$

a normal density with mean $\bar{y}_s$ and variance $\sigma^2/n$. Hence, $T = \bar{y}_s$ and $C = \bar{y}_s$.

EXAMPLE 2. Suppose that conditional on $\theta$, $Y$ is Bernoulli with success probability $\gamma = e^\theta/(1+e^\theta)$. The numerator of (2.4) becomes

$$e^{t\bar{y}_s n}/(1+e^t)^n = [e^t/(1+e^t)]^{n\bar{y}_s}[1-e^t/(1+e^t)]^{n-n\bar{y}_s}.$$

Hence, $V_n := e^{T_n}/(1+e^{T_n})$ has a Beta distribution converging to $\bar{y}_s$. Therefore, $C = \psi'(T) = \bar{y}_s$.

EXAMPLE 3. Suppose that conditional on $\theta$, $Y$ is Poisson with success rate $\lambda = e^\theta$. The numerator of (2.4) becomes $\exp[-n(e^t - t\bar{y}_s)]$. Hence, $V_n := e^{T_n}$ has the Gamma density $e^{-nv}v^{n\bar{y}_s-1}I(v>0)/\Gamma(n\bar{y}_s)$. Since $E(V_n) = \bar{y}_s$ and $\mathrm{Var}(V_n) \to 0$, $V_n$ converges to $\bar{y}_s$. Therefore, $C = \psi'(T) = \bar{y}_s$.

VERIFICATION OF (3.2). First, integrating (3.1) with respect to $\beta$, we obtain

$$\begin{aligned}
\pi(\boldsymbol{\theta}, r | y_s) \propto & \prod_{i=1}^{m} \prod_{j=1}^{n_i} \exp[\phi^{-1}\{y_{ij}\theta_i - \psi(\theta_i)\}] \\
& \times \prod_{i=1}^{m}\{h'(\theta_i)\} r^{m/2} \times \exp[-ar/2] r^{(b/2)-1} \\
& \times \int \exp\left[-(r/2)\sum_{i=1}^{m}\{h(\theta_i)-\beta\}^2\right] d\beta.
\end{aligned} \quad \text{(A.3)}$$

Then writing $\sum_{i=1}^{m}\{h(\theta_i)-\beta\}^2 = \sum_{i=1}^{m}\{h(\theta_i)-\bar{h}\}^2 + m(\beta-\bar{h})^2$, the last integral in (A.3) is proportional to $r^{-1/2}\exp[-(r/2)\sum_{i=1}^{m}\{h(\theta_i)-\bar{h}\}^2]$. Hence,

$$\begin{aligned}
\pi(\boldsymbol{\theta}, r | y_s) \propto & \prod_{i=1}^{m} \prod_{j=1}^{n_i} \exp[\phi^{-1}\{y_{ij}\theta_i - \psi(\theta_i)\}] \times \prod_{i=1}^{m}\{h'(\theta_i)\} \\
& \times r^{\{(m+b-1)/2\}-1} \exp\left\{-(r/2)\left[a + \sum_{i=1}^{m}\{h(\theta_i)-\bar{h}\}^2\right]\right\}.
\end{aligned} \quad \text{(A.4)}$$

Integrating (A.4) with respect to $r$ on $(0, \infty)$ and using the formula for the gamma integral,

$$\pi(\boldsymbol{\theta}|y_s) \propto \prod_{i=1}^{m} \prod_{j=1}^{n_i} \exp[\phi^{-1}\{y_{ij}\theta_i - \psi(\theta_i)\}] \times \prod_{i=1}^{m} h'(\theta_i)$$

(A.5)
$$\times \left[a + \sum_{i=1}^{m}\{h(\theta_i) - \bar{h}\}^2\right]^{-(b+m-1)/2}.$$

Integrating (A.5) with respect to $\theta_1, \ldots, \theta_{m-1}$, we obtain

(A.6) $\pi(\theta_m|y_s) \propto \prod_{j=1}^{n_m} \exp[\phi^{-1}\{y_{mj}\theta_m - \psi(\theta_m)\}] \times h'(\theta_m) \times g(\theta_m),$

where

(A.7)
$$g(\theta_m) = \int \prod_{i=1}^{m-1}\left\{h'(\theta_i)\prod_{j=1}^{n_i}\exp[\phi^{-1}\{y_{ij}\theta_i - \psi(\theta_i)\}]\right\}$$
$$\times \left[a + \sum_{i=1}^{m}\{h(\theta_i) - \bar{h}\}^2\right]^{-(b+m-1)/2} \prod_{i=1}^{m-1} d\theta_i.$$

Note that the function $g$ does not involve $n_m$ (which will be allowed to go to infinity). Then, from (A.6),

$$E[\psi'(\theta_m)|y_s] = \int \psi'(t)\pi(t|y_s)\,dt$$
$$= E_{T_{n_m}}[\psi'(T_{n_m})h'(T_{n_m})g(T_{n_m})]/E_{T_{n_m}}[h'(T_{n_m})g(T_{n_m})].$$

VERIFICATION OF (3.5). Note that

(A.8) $E[\psi'(\theta_m)|y_s] = \int \theta v_{n_m}(\theta|y_s)\,d\theta \Big/ \int v_{n_m}(\theta|y_s)\,d\theta,$

where

$v_{n_m}(\theta_m|y_s)$
$$= \int \left\{\sum_{i=1}^{m}\sum_{j=1}^{n_i}(y_{ij} - \theta_i)^2\right\}^{-(n_T/2+1)}\left\{a + \sum_{i=1}^{m}(\theta_i - \bar{\theta})^2\right\}^{-(b+m-1)/2} d\prod_{i=1}^{m-1}\theta_i$$
$$= \left\{n_m(\theta_m - \bar{y}_m)^2 + \sum_{j=1}^{n_m}(y_{mj} - \bar{y}_m)^2\right\}^{-(n_m/2+1)}$$
$$\times \int \left[\left\{n_m(\theta_m - \bar{y}_m)^2 + \sum_{j=1}^{n_m}(y_{mj} - \bar{y}_m)^2\right\}^{(n_m/2+1)}\right.$$
$$\times \left.\left(\left\{\sum_{i=1}^{m}\sum_{j=1}^{n_i}(y_{ij} - \theta_i)^2\right\}^{(n_T/2+1)}\right.\right.$$

$$\times \left\{ a + \sum_{i=1}^{m}(\theta_i - \bar{\theta})^2 \right\}^{(b+m-1)/2} \Bigg)^{-1} \Bigg] d\prod_{i=1}^{m-1} \theta_i$$

$$= \left\{ n_m(\theta_m - \bar{y}_m)^2 + \sum_{j=1}^{n_m}(y_{mj} - \bar{y}_m)^2 \right\}^{-(n_m/2+1)} \times l(\theta_m, y_s), \qquad \text{say.}$$

Write $h_{n_m}(\theta|y_s) = k(y_s)\{n_m(\theta_m - \bar{y}_m)^2 + \sum_{j=1}^{n_m}(y_{mj} - \bar{y}_m)^2\}^{-(n_m/2+1)}$, where $k(y_s)$ is such that $h_{n_m}$ is a density. Then, from (A.8),

$$(A.9) \quad E[\psi'(\theta_m)|y_s] = \int \theta l(\theta, y_s) h_{n_m}(\theta|y_s)\, d\theta \Big/ \int l(\theta, y_s) h_{n_m}(\theta|y_s)\, d\theta.$$

Next, we show that $h_{n_m}$ is the p.d.f. of a random variable that converges in probability to $\bar{y}_m$ as $n_m$ tends to infinity. Then the result that the hierarchical Bayes estimator converges to $\bar{y}_m$ follows by noting the boundedness of the function $\theta \Rightarrow \theta l(\theta, y_s)$.

For simplicity, write $n$ for $n_m$. Note that $h_n$ is a density of the form $c/\{n(x-\mu)^2 + d\}^{n/2+1}$. Clearly, the mean is $\mu$. To prove convergence in probability, we next show that for any $\varepsilon_1 > 0$,

$$(A.10) \quad \int_{\varepsilon_1}^{\infty} (nx^2 + d)^{-(n/2+1)}\, dx \Big/ \int_0^{\infty} (nx^2 + d)^{-(n/2+1)}\, dx \to 0.$$

Substituting $x = (d/n)^{1/2} y$, the above becomes equals $\int_{n^{1/2}\varepsilon}^{\infty}(y^2+1)^{-(n/2+1)} dy / \int_0^{\infty}(y^2+1)^{-(n/2+1)}\, dy$, where $\varepsilon := \varepsilon_1/d^{1/2}$. The numerator is bounded above by $\int_{n^{1/2}\varepsilon}^{\infty}(n^{1/2}\varepsilon y + 1)^{-(n/2+1)}\, dy$, which is $O(1/\{n^{3/2}(n^{1/2}\varepsilon)^n\})$. The denominator is bounded below by $\int (2e^y)^{-n}\, dy = O(1/n2^n)$. Hence, (A.10) follows.

**Acknowledgments.** We thank the Editor, an anonymous Associate Editor and two referees for their valuable comments which have led to significant improvement of an earlier version of the paper. We also record our deep appreciation to Ms. Huilin Li for her help with the computations.

Joint Program in Survey Methodology  
1218 Lefrak Hall  
University of Maryland at College Park  
College Park, Maryland 20742  
USA  
E-mail: plahiri@survey.umd.edu

Department of Mathematical Sciences  
University of Liverpool  
Liverpool L69 7ZL  
United Kingdom  
E-mail: k.mukherjee@liverpool.ac.uk